\documentclass[11pt]{article}

\usepackage{fullwidth}
\usepackage{xypic}
\usepackage{xy}
\usepackage[top=1.5in, bottom=1.5in, left=1in, right=1in]{geometry}
\usepackage{amsgen}
\usepackage{amsmath}
\usepackage{amstext}
\usepackage{amsbsy}
\usepackage{amsopn}
\usepackage{amsfonts}
\usepackage{amssymb}
\usepackage{eepic}
\usepackage{graphicx}
\usepackage{epsf}
\usepackage{pstricks}
\usepackage{xcolor}

\xyoption{all}

\def\Box{\square}

\def\tra#1{\smash{\mathop{\mid\kern
-1pt\joinrel\relbar\joinrel\relbar}\limits^{*}_{#1}}}
\def\longtra#1{\smash{\mathop{\mid\kern
-1pt\joinrel\relbar\joinrel\relbar\joinrel\relbar}\limits^{*}_{#1}}}
\def\vlongtra#1{\smash{\mathop{\mid\kern
-1pt\joinrel\relbar\joinrel\relbar\joinrel\relbar\joinrel\relbar}\limits^{*}_{#1}}}
\def\vvlongtra#1{\smash{\mathop{\mid\kern
-1pt\joinrel\relbar\joinrel\relbar\joinrel\relbar\joinrel\relbar\joinrel\relbar}\limits^{*}_{#1}}}
\def\vvvlongtra#1{\smash{\mathop{\mid\kern
-1pt\joinrel\relbar\joinrel\relbar\joinrel\relbar\joinrel\relbar\joinrel\relbar\joinrel\relbar}\limits^{*}_{#1}}}
\def\etra#1{\smash{\mathop{\mid\kern
-1pt\joinrel\relbar\joinrel\relbar}\limits_{#1}}}

\def\A{{\cal{A}}}

\def\iff{\Leftrightarrow}
\def\Rw{\Rightarrow}

\def\wt{\widetilde}

\def\N{\mathbb{N}}

\def\mt{\mbox{MT}}

\def\P{{\cal{P}}}

\def\R{{\mathcal{R}}}

\def\Z{\mathbb{Z}}

\def\p{\varphi}

\def\inv{^{-1}}
\def\la{\langle}
\def\ra{\rangle}

\def\bi{\begin{itemize}}
\def\ei{\end{itemize}}
\def\beq{\begin{equation}}
\def\eeq{\end{equation}}

\def\J{{\mathcal J}}

\def\xr{\xrightarrow}

\newtheorem{T}{Theorem}[section]
\newcommand{\bt}{\begin{T}}
\newcommand{\et}{\end{T}}
\newcommand{\ftd}{$\square$\end{T}}

\newtheorem{Proposition}[T]{Proposition}
\newcommand{\bp}{\begin{Proposition}}
\newcommand{\ep}{\end{Proposition}}
\newcommand{\fpd}{$\square$\end{Proposition}}

\newtheorem{Lemma}[T]{Lemma}
\newcommand{\bl}{\begin{Lemma}}
\newcommand{\el}{\end{Lemma}}
\newcommand{\fld}{$\square$\end{Lemma}}

\newtheorem{Corol}[T]{Corollary}
\newcommand{\bc}{\begin{Corol}}
\newcommand{\ec}{\end{Corol}}
\newcommand{\fcd}{$\square$\end{Corol}}

\newtheorem{Result}[T]{Result}
\newcommand{\br}{\begin{Result}}
\newcommand{\er}{\end{Result}}
\newcommand{\frd}{$\square$\end{Result}}

\newtheorem{Example}[T]{Example}
\newcommand{\be}{\begin{Example}}
\newcommand{\ee}{\end{Example}}

\newtheorem{Problem}[T]{Problem}
\newcommand{\bq}{\begin{Problem}}
\newcommand{\eq}{\end{Problem}}

\newtheorem{Remark}[T]{Remark}
\newcommand{\brem}{\begin{Remark}}
\newcommand{\erem}{\end{Remark}}

\newcommand{\proof}
   {\par\medbreak\noindent{\bf Proof}.\enspace}

\newcommand{\qed}{
$\Box$
\par\bigbreak}

\title{Finitely presented inverse semigroups with finitely many idempotents in each $\mathcal D$-class and non-Hausdorff universal groupoids}
\author{{\bf Pedro V. Silva}\\ $ $\\
{\em Centro de
Matem\'{a}tica, Faculdade de Ci\^{e}ncias, Universidade do
Porto,}\\ {\em R. Campo Alegre 687, 4169-007 Porto, Portugal}\\
{\em email:} pvsilva@fc.up.pt\\
$ $\\
{\bf Benjamin Steinberg}\\ $ $\\
{\em Department of Mathematics, City College of New York,}\\
{\em Convent Avenue at 138th Street, New York, New York 10031, USA}\\
{\em email:} bsteinberg@ccny.cuny.edu}

\date{\today}

\begin{document}
\maketitle

\begin{center}\small

\bigskip

\end{center}

\abstract{The complex algebra of an inverse semigroup with finitely many idempotents in each $\mathcal D$-class is stably finite by a result of Munn.  This can be proved fairly easily using $C^*$-algebras for inverse semigroups satisfying this condition that have a Hausdorff universal groupoid, or more generally for direct limits of inverse semigroups satisfying this condition and having Hausdorff universal groupoids.  It is not difficult to see that a finitely presented inverse semigroup with a non-Hausdorff universal groupoid cannot be a direct limit of inverse semigroups with Hausdorff universal groupoids.  We construct here countably many non-isomorphic finitely presented inverse semigroups with finitely many idempotents in each $\mathcal D$-class and non-Hausdorff universal groupoids.  At this time there is not a clear $C^*$-algebraic technique to prove these inverse semigroups have stably finite complex algebras.
}

\section{Introduction}

Given an inverse semigroup $S$, we denote by $E(S)$ its semilattice of idempotents and by $\leq$ its natural partial order (hence $a \leq b$ if and only if $a = eb$ for some $e \in E(S)$). Given $s \in S$, let
$\lambda_S(s) = \{ t \in S \mid t \leq s\}$ and let
$\mu_S(s)$ denote the set of maximal elements of $\lambda_S(s) \cap E(S)$ for $\leq$.

Paterson~\cite{Paterson} associated an \'etale groupoid $\mathcal G(S)$ to every inverse semigroup $S$, called its universal groupoid, and showed that the $C^*$-algebra of the inverse semigroup is isomorphic to the $C^*$-algebra of its universal groupoid.  The second author later generalized Paterson's result by showing that if $K$ is any commutative ring with unit, then the semigroup algebra $KS$ of $S$ is isomorphic to a certain convolution algebra of $K$-valued functions on the groupoid $\mathcal G(S)$~\cite{groupoidalgebra}.  There is now a well-developed theory of \'etale groupoid algebras which has proven useful for studying inverse semigroup algebras.  Both the algebra and $C^*$-algebra of an \'etale groupoid are best behaved when the groupoid is Hausdorff. In particular, the theory of traces on groupoid $C^*$-algebras seems to only be well developed in the Hausdorff case where there is a faithful conditional expectation mapping to the algebra of continuous functions vanishing at infinity on the unit space.   Universal groupoids of inverse semigroups are not always Hausdorff.  It was shown in~\cite{groupoidalgebra} that $\mathcal G(S)$ is Hausdorff if and only if $\lambda_S(s)\cap E(S)$ is finitely generated as an ideal in $E(S)$, that is, there is a finite set $F\subseteq \lambda_S(s)\cap E(S)$ such that if $e\in \lambda_S(s)\cap E(S)$, then $e\leq f$ for some $f\in F$.  In particular, if $\mathcal G(S)$ is Hausdorff, then $\mu_S(s)$ is finite.  For example, if $G$ is a nontrivial group and $E$ is the semilattice consisting of a zero element and a countably infinite set of orthogonal idempotents, then $S=G\cup E$ is an inverse monoid, where $G$ acts trivially on the left and right of $E$, with a non-Hausdorff universal groupoid.  Indeed, if $1\neq g\in G$, then $\lambda_S(g)\cap E(S)=E$ which has infinitely many maximal elements.  Note that $S$ is a Clifford inverse monoid: each $\mathcal D$-class of $S$ contains a single idempotent.  But $S$ is not finitely generated and each of its finitely generated inverse subsemigroups does have a Hausdorff universal groupoid and hence $S$ is a direct limit of Clifford monoids with Hausdorff universal groupoids.

The second author has recently initiated a study of stable finiteness of \'etale groupoid algebras~\cite{traces} and, in particular, recovered a result of Munn showing that if $S$ is an inverse semigroup whose $\mathcal D$-classes have finitely many idempotents, then $KS$ is stably finite for any field $K$ of characteristic $0$~\cite{MunnDirect}.  Recall that a ring $R$ is stably finite if $M_n(R)$ does not contain a copy of the bicyclic monoid as a subsemigroup for any $n\geq 1$.    In the case that $S$ has a Hausdorff universal groupoid, this can be deduced using the theory of $C^*$-algebras, but in the non-Hausdorff case one needs to work around this.   However, since the stable finiteness result can be reduced to the case of finitely generated inverse semigroups with finitely many idempotents in each $\mathcal D$-class, it becomes of interest to know whether there are examples of such finitely generated inverse semigroups with non-Hausdorff universal groupoids.  The easiest way to guarantee that each $\mathcal D$-class of $S$ contains finitely many idempotents is to impose the stronger condition that each $\R$-class of $S$ is finite (since each idempotent $f$ in the $\mathcal D$-class of $e$ is of the form $s^{-1}s$ with $s$ in the $\R$-class of $e$).

This paper was then  motivated by the following question:

\bq
\label{question}
Is there a finitely generated  semigroup $S$ such that:
\bi
\item[(i)] every $\R$-class of $S$ is finite;
\item[(ii)] $\mu_S(s)$ is infinite for some $s \in S$?
\ei
\eq

In Section~\ref{sg} we present a construction which provides uncountably many nonisomorphic 3-generated inverse monoids satisfying the conditions of Problem~\ref{question}.

In Section~\ref{sp} we present a construction which provides infinitely many nonisomorphic finitely presented 3-generated inverse monoids satisfying the conditions of Problem~\ref{question}.

It is not difficult to show (see~\cite{traces}) that if $S$ is finitely presented and has a non-Hausdorff universal groupoid, then $S$ cannot be written as a direct limit of inverse semigroups with Hausdorff universal groupoids.  Since stable finiteness is preserved under direct limits, to really show that we cannot reduce stable finiteness of $KS$ to the case that $S$ has a Hausdorff universal groupoid, it is important to have a finitely presented inverse semigroup satisfying the conditions of Problem~\ref{question}.

\section{Preliminaries}

The reader is assumed to have basic knowledge of inverse semigroup theory and automata theory, being respectively referred to~\cite{Law} and~\cite{Sak} for that purpose. Since the inverse semigroups we construct are actually inverse monoids, all the relevant definitions are presented in the monoid version.

\subsection{Inverse automata}

Given a finite alphabet $A$, we
denote by $A\inv$ a set of formal inverses of $A$ and write $\wt{A} = A\cup
A\inv$.
An {\em inverse automaton} over the alphabet $\widetilde{A}$ is a
structure of the form ${\cal{A}} = (Q,i,t,E)$, where
\bi
\item
$Q$ is the set of vertices,
\item
$i, t \in Q$ are the initial and terminal vertices, respectively,
\item
$E \subseteq Q \times \wt{A} \times Q$ is the set of edges,
\ei
satisfying the following properties:
\bi
\item
{\em deterministic}: $(p,a,q),(p,a,q') \in E \Rw q = q'$;
\item
{\em involutive}: $(p,a,q) \in E \iff (q,a^{-1},p) \in E$;
\item
{\em trim}: every vertex lies in some path from $i$ to $t$.
\ei

If $i = t$, we refer to it as the {\em basepoint} of $\A$. If we do not specify the initial and terminal vertices, we have an {\em inverse graph}.

Assume now that $\A$ is involutive and trim, but not deterministic. A {\em folding} operation on $\A$ consists on identifying two distinct edges of the form $p \xleftarrow{a} q \xr{a} r$ (identifying also the inverse edges $p \xr{a\inv} q \xleftarrow{a\inv} r$). If $\A$ is finite and we fold enough edges, we end up obtaining a finite inverse automaton $\A'$, which we say is obtained from $\A$ by {\em complete folding}.

Is this operation confluent? That is, does the inverse automaton depend on the sequence of foldings? Given $u,v \in \wt{A}^*$, write $u \xr{*} v$ if $v$ can be obtained from $u$ by successively erasing factors of the form $aa\inv$ $(a \in \wt{A})$. It is easy to check that
$$L(\A') = \{ v \in \wt{A}^* \mid u \xr{*} v\mbox{ for some } u \in L(\A)\}.$$
Hence the language of the inverse automaton $\A'$ is completely determined by $\A$. Since inverse automaton are known to be minimal (see e.g.~\cite{BS}), it follows that $\A'$ is itself determined by $\A$ and so the folding process is confluent.

A {\em Dyck word} on $\wt{A}$ is some $w \in \wt{A}^*$ satisfying $w \xr{*} 1$. These are the words representing the identity in the free group on $A$, and play also an important role in the theory of inverse semigroups as we soon shall see.

\subsection{Free inverse monoids}

Let $A$ be a nonempty alphabet. We extend $^{-1}\colon A \to A^{-1}: a \mapsto a^{-1}$ to an involution on the free monoid
$\widetilde{A}^*$ through
$$1\inv = 1,\quad (a^{-1})^{-1} = a,\quad (uv)^{-1} = v^{-1}u^{-1}\hspace{1cm} (a \in A;\;
u,v \in \widetilde{A}^+)\, .$$
The {\em free inverse monoid on} $A$ is the quotient $FIM_A = \wt{A}^*/\rho$,
where $\rho$ is the congruence on $\wt{A}^*$ generated by the relation
$$\{ (ww\inv w,w) \mid w \in \wt{A}^*\} \cup
\{ (uu\inv vv\inv,vv\inv uu\inv) \mid u,v \in \wt{A}^*\}.$$
known as the {\em Wagner congruence} on $\wt{A}^*$.

W. D. Munn provided in~\cite{Mun} an elegant normal form for $FIM_A$ using inverse automata (see also~\cite{Sch} by Scheiblich).

Given $w = a_1\ldots a_n \in \wt{A}^*$ $(a_i \in \wt{A})$, let $\mathrm{Lin}(w)$ denote the {\em linear automaton} of $w$:
$$\xymatrix{
\ar[r] & q_0 \ar[r]^{a_1} & q_1 \ar[r]^{a_2} & \ldots \ar[r]^{a_n} & q_n \ar[r] &
}$$
which contains also the inverse edges (to make it involutive). The {\em Munn tree}
of $w$ is the finite inverse automaton
$\mt(w)$ obtained by completely folding $\mathrm{Lin}(w)$. This provides the following solution for the word problem of $FIM_A$:

\bt
\label{wpfree}
{\rm~\cite{Mun}}
For all $u,v \in \wt{A}^*$, the following conditions are equivalent:
\bi
\item[(i)] $u\rho = v\rho$;
\item[(ii)] ${\rm MT}(u) \cong {\rm MT}(v)$;
\item[(iii)] $u \in L({\rm MT}(v))$ and
  $v \in L({\rm MT}(u))$.
\ei
\et

Such Munn trees are precisely those finite inverse automata on $\wt{A}$ whose underlying undirected graph is a tree (when we consider only the edges labeled by $A$).

It is easy to see that, given $w \in
\wt{A}^*$, we have
$$w\rho \in E(FIM_A) \iff \mbox{ $\mathrm{Lin}(w)$ has a basepoint $\iff w$ is a Dyck word}.$$

\subsection{Inverse monoid presentations}

A (finite) inverse monoid presentation is a formal
expression of the form
$\P = \langle A \mid R\rangle$, where $A$ is a (finite) alphabet and $R$ is
a (finite) subset of $\wt{A}^* \times \wt{A}^*$. We usually describe the relations in $R$ as formal equalities $r = s$.

Let $\tau = (\rho
\cup R)^{\sharp}$ be the congruence on $\wt{A}^*$ generated by the relation $\rho \cup R$. The quotient
$S = \wt{A}^*/\tau$ is the inverse  monoid defined by
$\P$. Equivalently, we might write $S = FIM_A/R^{\sharp}$, viewing $R$ as a relation on $FIM_A$. We denote by $\p\colon \wt{A}^* \to S$ the canonical homomorphism.

Stephen devised in~\cite{Ste} an approach which is the most useful tool known to date to deal with inverse  monoid presentations. The {\em Cayley graph} of $S$ with respect to the generating set $A$, denoted by Cay$_A(S)$, has vertex set $S$ and edges of the form $s \xr{a} s(a\p)$ for all $s \in S$ and $a \in \wt{A}$. This is not in general an involutive graph (it would be if $S$ is a group). But the strongly connected components of Cay$_A(S)$ are actually inverse automata. And these correspond to the various $\R$-classes of $S$ (so there exist paths
$$\xymatrix{
s_1 \ar@/^/[rr]^u && s_2  \ar@/^/[ll]^v
}$$
in Cay$_A(S)$ for some $u,v \in \wt{A}^*$ if and only if $s_1s_1\inv = s_2s_2\inv$).
The {\em Sch\"{u}tzenberger graph} of $w \in \wt{A}^*$, denoted by $S\Gamma(w)$,  is the strongly connected component of Cay$_A(S)$ containing $w\p$ (that is, the induced subgraph having the $\R$-class of $w\p$ as set of vertices). The {\em Sch\"{u}tzenberger automaton} of $w \in \wt{A}^*$, denoted by $\A(w)$, is obtained from $S\Gamma(w)$ by setting $(ww\inv)\p$ as initial vertex and $w\p$ as terminal vertex.
We may also write $S\Gamma(u\p) = S\Gamma(u)$ or $\A(u\p) = \A(u)$ if it suits us.

Stephen proved the following theorems, which generalize those of Munn:

\bt
\label{sanpo}
{\rm~\cite{Ste}}
For all $u,v \in \wt{A}^*$, the following conditions are equivalent:
\bi
\item[(i)] $u\p \geq v\p$;
\item[(ii)] $u \in L(\A(v))$.
\ei
\et

\bt
\label{sawp}
{\rm~\cite{Ste}}
For all $u,v \in \wt{A}^*$, the following conditions are equivalent:
\bi
\item[(i)] $u\p = v\p$;
\item[(ii)] $\A(u) \cong \A(v)$;
\item[(iii)] $u \in L(\A(v))$ and
  $v \in L(\A(u))$.
\ei
\et

Since $u\p \in E(S)$ if and only if $(uu\inv)\p = u\p$, it follows that the idempotents of $S$ are characterized by having a basepoint at their Sch\"{u}tzenberger automaton.

It follows from Theorem~\ref{sawp} that the word problem is decidable for $\P$ if membership is decidable in the languages of its Sch\"{u}tzenberger automata. In general, it is not. But Stephen devised a procedure which brings positive results in many important cases. We describe it now.

Suppose that $\A$ is a finite involutive automaton and $r = s$ is a relation in $R$ such that there exists in $\A$ a path of the form $p \xr{r} q$ but no path $p \xr{s} q$ (or vice-versa). If we glue to $\A$ a ``path" $p \xr{s} q$ with the corresponding inverse edges (to keep it involutive), we say that this new automaton $\A'$ is obtained from $\A$ by an {\em expansion} (expanding through the relation $r = s$).

If $L(\A) \subseteq L(\A(w))$ and $\A'$ is obtained from $\A$ by either folding or through an expansion through some relation in $R$, then $L(\A') \subseteq L(\A(w))$.
We call any finite inverse automaton $\A$ such that $w \in L(\A) \subseteq L(\A(w))$ an {\em approximate automaton} of $w$ (for the presentation $\P$). This is the case of any finite inverse automaton
obtained from $\mt(w)$ through a finite sequence of foldings and expansions. If $\A$ is an approximate automaton of $w$ which admits neither foldings nor expansions, then $\A \cong \A(w)$.

Assume now that the $\R$-classes of $S$ are all finite. Then the Sch\"{u}tzenberger automata are all finite as well. Can we compute then? If $R$ is finite, the answer is given through Stephen's sequence, which provides a systematic way of building approximate automata. Take $\A_1(w) = \mt(w)$. If $\A_n(w)$ is defined, let $\A'_{n+1}(w)$ be obtained by performing {\em simultaneously} all the possible expansions of $\A_n(w)$. Since $R$ is finite, $\A'_{n+1}(w)$ is a finite involutive automaton. Now define $\A_{n+1}(w)$ by completely folding $\A'_{n+1}(w)$. We get then a sequence $(\A_n(w))_n$ of approximate automata of $w$. If $\A(w)$ is finite, it will eventually show up as a member of the sequence.

This is still true for infinite $R$ if:
\bi
\item
we are sure of $S$ having only finite $\R$-classes;
\item
only finitely many expansions can be applied to each $\A_n(w)$.
\ei

Anyway, in the general case of an arbitrary presentation, $v \in L(\A_T(u))$ implies that $v$ is recognized by some approximate automaton of $u$; in fact given any approximate automaton of $u$, one can perform a finite sequence of foldings and expansions to obtain an approximate automaton recognizing $v$.

The results in this subsection will be used throughout the paper without further reference. The reader is referred to~\cite{Ste} for more details on Sch\"{u}tzenberger automata and Stephen's sequences.

\section{Finitely generated examples}
\label{sg}

In this section we present a construction which provides uncountably many nonisomorphic 3-generated inverse monoids satisfying the conditions of Problem~\ref{question}.

We start by proving the following:

\bp
\label{con}
Let $S$ be an inverse monoid where every $\R$-class is finite.
Let $\la A \mid R\ra$ be an inverse monoid presentation of $S$. Let $E \subseteq E(S)\setminus \{ 1\}$ and let $w_e \in \wt{A}^+$ be a Dyck word representing
$e$ for every $e \in E$. Let $T$ be the inverse semigroup defined by the inverse monoid presentation
\beq
\label{con1}
\la A \cup \{ b\} \mid R,\, w_eb = w_e\, (e \in E)\ra,
\eeq
where $b$ is a new letter. Then:
\bi
\item[(i)]
every $\R$-class of $T$ is finite;
\item[(ii)]
if no two distinct elements of $E$ are $\J$-comparable in $S$ and $E$ is infinite, then $\mu_T(b)$ is infinite.
\ei
\ep

\proof
(i) Note that an idempotent $e$ can always be represented by some Dyck word since $e = ee\inv$.
Write $B = A \cup \{ b\}$ and let $\psi\colon \wt{A}^* \to S$ and $\p\colon \wt{B}^* \to T$ be the canonical homomorphisms.

We shall use the notation $\A_S(u)$ and $\A_T(u)$ to denote the Sch\"utzenberger automaton of $u$ relative to the presentation of $S$ and $T$, respectively.

Let $u \in \wt{B}^*$. Let $Q'$ be the set of vertices of $\mt(u)$ which admit a $b$-loop (at their image) in $\A_T(u)$. Let $\A_1$ be the finite inverse automaton obtained from $\mt(u)$ by adjoining a $b$-loop at each $q \in Q'$ followed by complete folding. Then $\A_1$ is an approximate automaton of $u$ with the property that each $b$-edge that is not a loop comes from $\mt(u)$.
If we remove the $b$-edges from $\A_1$, then we have essentially a disjoint union of Munn trees with edges in $\wt{A}$ since we only folded $b$-edges. Since every $\R$-class of $S$ is finite, we can turn each of these Munn trees into a finite inverse automaton admitting no $R$-expansions by a finite series of $R$-expansions and foldings. The resulting automaton
$\A_2$ is a finite inverse automaton admitting no $R$-expansions, and is still an approximate automaton of $u$.
Now let $\A_3$ be obtained from $\A_2$ by adjoining a $b$-loop at each vertex $p$ admitting a path $p \xr{w_e\inv} \ldots$ for some $e \in E$ (if it does not exist already). Adjoining these new $b$-loops does not allow any folding: if $p \xr{b^{\varepsilon}} q$ would be an edge of $\A_2$ for $\varepsilon = \pm 1$, then we have necessarily $p \in Q'$ because $p$ is doomed to host a $b$-loop in $\A_T(u)$, and then we already have a $b$-loop at $p$ in $\A_1$. Hence $\A_3$ admits no folding and it certainly admits no expansions of any sort. Since $\A_3$ is an approximate automaton of $u$, it follows that $\A_3 = \A_T(u)$, which is therefore finite.

Therefore $\A_T(u)$ must be finite in all cases and so is $\R_{u\p}$.

(ii) We start by noting that the homomorphism $\eta\colon  S \to T$ extending the identity mapping on $A$ is an embedding. It is immediate that we have a homomorphism $\theta\colon T \to S$ defined by
$$b\theta = 1,\quad a\theta = a \; (a \in A).$$
Then $\theta\eta=1_S$  (as it is the identity on $A$), and so $\eta$ is an embedding of $S$ into $T$.

Let $e \in E$. To construct $\A_T(w_e)$, we start by turning $\mt(w_e)$ into $\A_S(w_e)$, which is an approximate automaton of $w_e$ for the presentation (\ref{con1}). Expanding through the relation $w_eb = w_e$ and subsequent folding will produce a $b$-loop at the basepoint $i$ and possibly at other vertices where a path labeled by $w_e$ (necessarily a loop) can be read.

Suppose now that we have some path of the form $p \xr{w_{e'}} q$ in $\A_S(w_e)$ with $e' \in E \setminus \{ e\}$. Then $uw_{e'}v \in L(\A_S(w_e))$ for some $u,v \in \wt{A}^*$ and so $(u\psi)e'(v\psi)e = e$, yielding $e' \geq_{\J} e$. This contradicts the assumption that no two distinct elements of $E$ are $\J$-comparable in $S$, hence $\A_S(w_e)$ admits no path of the form $p \xr{w_{e'}} q$ and so $\A_T(w_e)$ is indeed $\A_S(w_e)$ with a few $b$-loops attached (namely at the basepoint).

In view of the relation $w_eb = w_e$, we get $w_e\p \leq b\p$ and so $w_e\p \in \lambda_T(b\p) \cap E(T)$.
Suppose now that $w_e\p \leq f \in \lambda_T(b\p) \cap E(T)$. Let $v \in \wt{B}^*$ be a Dyck word representing $f$.

On the one hand, $w_e\p \leq v\p$ implies $L(\A_T(v)) \subseteq L(\A_T(w_e))$. On the other hand,
since $f \leq b\p$, we have $v\p =f = f(b\p) = (v\p)(b\p) =(vb)\p$.
Thus $\A_T(vb) \cong \A_T(v)$ and so there exists a $b$-loop at $i$ in $\A_T(vb)$.

It follows that at some point in the Stephen's sequence of $vb$, we must have had some expansion through some relation of the form $w_{e'}b = w_{e'}$. Hence $w_{e'}$ labels some path in $\A_T(v)$ and so $xw_{e'}y \in L(\A_T(v)) \subseteq L(\A_T(w_e))$
for some $x,y \in \wt{B^*}$.

We have seen above that there is no path in $\A_S(w_e)$ (and consequently neither in $\A_T(w_e)$) labeled by some $w_{e'}$ for $e' \in E \setminus \{ e\}$, hence $\A_T(v)$ must admit some loop labeled by $w_e$ at some vertex $q_1$.

Let $i = q_0 \xr{w} q_1$ be a path in $\A_T(v)$. For those familiar with the bicyclic monoid, if we put $x=(vww_e)\p$, then $xx\inv =f$ and $x\inv x = e(w\inv\p)f(w\p) e \leq e\leq f$.  So if $e<f$, then $x,x\inv$ generate a copy of the bicyclic monoid with identity $f$ and hence the $\R$-class of $f$ is infinite, contradicting part (i).  For those not familiar with the bicyclic monoid, here is a direct proof.
Then $ww_ew\inv \in L(\A_T(v)) \subseteq L(\A_T(w_e))$ and so $(ww_ew\inv)\p \geq w_e\p$.
Assume that there exists a path
$$\xymatrix{
q_0 \ar[r]^{w} & q_1 \ar[r]^{w} & \ldots \ar[r]^{w} & q_j \ar@(u,ur)^{w_e}
}$$
in $\A_T(v)$ for some $j \geq 1$. Then $w^jw_ew^{-j} \in L(\A_T(v))$ yields $(w^jw_ew^{-j})\p \geq v\p$ and so
$$(w^{j+1}w_ew^{-(j+1)})\p = (w^j(ww_ew\inv)w^{-j})\p \geq (w^jw_ew^{-j})\p \geq v\p.$$
Thus $w^{j+1}w_ew^{-(j+1)} \in L(\A_T(v))$.
By induction, it follows that we have a path of the above form in $\A_T(v)$ for every $j \geq 1$. Since $\A_T(v)$ is finite by part (i) and is an inverse automaton, we must have $q_j = q_0$ for some $j \geq 1$. But then $w_e$ labels a loop at $q_0 = i$ in $\A_T(v)$ and so $w_e\p \geq v\p = f$.

Thus $f = w_e\p$ and so $w_e\p \in \mu_T(b\p)$ for every $e \in E$.

Since $E$ is infinite and $\eta\colon S\to T$ induced by the identity on $A$ is an embedding, $\{ w_e\p \mid e \in E\}$ is an infinite subset of $\mu_T(b\p)$.
\qed

Now we can use Proposition~\ref{con} to produce examples which answer positively Problem~\ref{question}:

\be
\label{fgex}
Let $I \subseteq \N \setminus \{ 0 \}$ be infinite and let $T_I$ be defined by the inverse monoid presentation
\beq
\label{fgex1}
\la a,b,c \mid ab^ia = ab^iac\; (i \in I)\ra
\eeq
Then:
\bi
\item[(i)] all the $\R$-classes of $T_I$ are finite;
\item[(ii)] $\mu_{T_I}(c)$ is infinite.
\ei
\ee

Indeed, let $S = FIM_{\{ a,b\}}$, which 
has finite $\R$-classes~\cite{Mun}. Then $E_I = \{ ((ab^ia)\inv(ab^ia))\rho \mid i \in I\}$ is an infinite subset of $E(S)$. If $i,j \in I$ are distinct, then there is no path labeled by $ab^ia$ in $\mt(ab^ja)$ and vice-versa, so no two distinct elements of $E_I$ are $\J$-comparable in $S$.
Consider now the inverse monoid presentation
$$\la a,b,c \mid (ab^ia)\inv(ab^ia) = (ab^ia)\inv(ab^ia)c\; (i \in I)\ra.$$
In any inverse semigroup, $u = uc$ is equivalent to $u\inv u = u\inv uc$, hence this presentation is equivalent to (\ref{fgex1}).

Now it follows from Proposition~\ref{con} that all the $\R$-classes of $T_I$ are finite and $\mu_{T_I}(c)$ is infinite.

\medskip

We can now use Example~\ref{fgex} to prove the following:

\bp
\label{umfg}
There exist uncountably many nonisomorphic 3-generated inverse monoids satisfying the conditions of Problem~\ref{question}.
\ep

\proof
For each infinite $I \subseteq \N \setminus \{ 0 \}$, let $T_I$ be defined by the inverse monoid presentation
(\ref{fgex1}). We have shown in Example~\ref{fgex} that $T_I$ satisfies the conditions of Problem~\ref{question}. Since a countably infinite set contains uncountably many infinite subsets, it suffices to show that
$T_I \not\cong T_J$ for distinct $I,J \subseteq \N \setminus \{ 0 \}$.

Out of symmetry, we may assume that $i \in I\setminus J$. Suppose that $\theta\colon T_I \to T_J$ is an isomorphism. It is easy to see that each generating set of $T_J$ must contain
\bi
\item
$a$ or $a\inv$;
\item
$b$ or $b\inv$;
\item
$c$ or $c\inv$.
\ei
Hence any minimal generating set of $T_J$ is necessarily of the form $\{ a^{\varepsilon},b^{\delta}, c^{\gamma} \}$ with $\varepsilon, \delta, \gamma = \pm 1$. Thus we may assume that $\theta$ is induced by some bijection $\{ a,b,c\} \to \{ a^{\varepsilon},b^{\delta}, c^{\gamma} \}$.

Since $ab^ia = ab^iac$ is a relation of the presentation of $T_I$, then
$(ab^ia)\theta = (ab^iac)\theta$ holds in $T_J$. Hence $c\theta$ labels a loop in $\A_{T_J}((ab^ia)\theta)$ and so $c\theta =
c^{\gamma}$ necessarily.

Since $ab^ia = ab^iac$ holds in $T_I$, then $(ab^ia)\theta = (ab^ia)\theta c^{\gamma}$ must hold in $T_J$.
The only way of enabling an expansion in $\mt((ab^ia)\theta)$ is if $a\theta = a$, $b\theta = b$ and $i \in J$, a contradiction. Therefore $T_I \not\cong T_J$.
\qed

\section{Finitely presented examples}
\label{sp}

In this section we present a construction which provides infinitely many nonisomorphic finitely presented 3-generated inverse monoids satisfying the conditions of Problem~\ref{question}.

\be
\label{fpex}
For each $t \geq 2$, let $S_t$ be defined by the inverse monoid presentation
\beq
\label{fpex1}
\la a,b,c \mid ca = a, 
cb^{-t}c\inv b^{t} = cb^{-t}b^{t}c\inv\ra.
\eeq
Then:
\bi
\item[(i)] every $\R$-class of $S_t$ is finite;
\item[(ii)] $\mu_{S_t}(aca\inv)$ is infinite.
\ei
\ee

Let us check these facts.

(i) Write $S = S_t$ and $A = \{ a,b,c\}$. Let $\p\colon \wt{A}^* \to S$ be the canonical homomorphism.
For every $x \in A$, let $\pi_x\colon \wt{A}^* \to \Z$ be the homomorphism defined by
$$y\pi_x = \left\{
\begin{array}{ll}
1&\mbox{ if }y = x\\
-1&\mbox{ if }y = x\inv\\
0&\mbox{ if }y \in \wt{A} \setminus \{ x,x\inv\}
\end{array}
\right.$$

Let $u \in \wt{A}^*$. It suffices to show that $S\Gamma(u)$ possesses only finitely many edges labeled by $a,b,c$.

If we perform an expansion involving the letter $a$, the subsequent folding prevents the number of $a$-edges to increase. Hence the number of $a$-edges in $S\Gamma(u)$ is bounded by the number of $a$-edges in $\mt(u)$.

Dealing with the $b$-edges is harder because its number may increase through the Stephen's sequence of $u$. We start by proving a few remarks.

\beq
\label{fpex3}
\mbox{If $v$ labels a loop in $S\Gamma(u)$, then }v\pi_b = 0.
\eeq

Inspection of the defining relations shows that $\pi_b$ induces a homomorphism $\overline \pi_b\colon S_t\to \mathbb Z$.  If $v$ reads a loop in $S\Gamma(u)$, then $s(v\p) =s$ for some $s\in S_t$, and so $s\overline \pi_b +v\pi_b = s\overline \pi_b$, whence $v\pi_b=0$.

Let $v \in \wt{A}^+$. We call a vertex (respectively edge) of $S\Gamma(v)$ {\em original} if it corresponds to some vertex (respectively edge) of $\mt(v)$. We show that:

\beq
\label{fpex4}
\mbox{Every vertex $p$ of $S\Gamma(u)$ admits some path $p \xr{b^m} q$ with $m \geq 0$ and $q$ original}.
\eeq

Once again, this property holds trivially for $\mt(u)$ and is preserved through folding. What about expansions?
If we expand through the relation $ca = a$ and fold once to make the new $c$ label a loop, the property still holds.
Hence we may assume that we are expanding through the relation $cb^{-t}c\inv b^{t} = cb^{-t}b^{t}c\inv $. The only way of getting new vertices is described by the following picture (where some appropriate folding is also assumed to limit the appearance of new vertices to a minimum):
$$
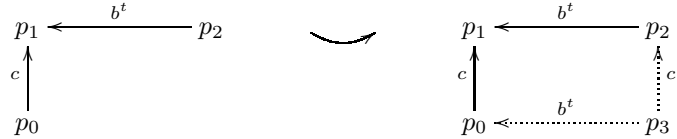
\xymatrix{
p_1  && p_2 \ar[ll]_{b^t} &  \ar@/_/[r] & & p_1 && p_2 \ar[ll]_{b^t}\\
p_0 \ar[u]^c &&&&&p_0 \ar[u]^c &&p_3  \ar@{.>}[ll]_{b^t} \ar@{.>}[u]_c
}$$
If the property holds for the vertex $p_0$, it must also hold for the new vertices ($p_4$ and $p_5$ in the picture), thus the property
is preserved through 
both sorts of expansions. It follows that the property holds throughout the whole Stephen's sequence of $u$. Since
every vertex $p$ of $S\Gamma(u)$ must originate from some term of the Stephen's sequence, then (\ref{fpex4}) holds.

Let $X$ be the set of $b$-edges of $S\Gamma(u)$. Let $\sigma$ be the equivalence relation on $X$ generated by relating edges
$$p \xr{b} q \quad \mbox{and} \quad p' \xr{b} q'$$
whenever $S\Gamma(u)$ admits a path
$$\xymatrix{
p \ar[r]^{b^s} & \bullet \ar[r]^c & \bullet \ar[r]^{b^{-s}} & p'
}$$
for some $s \in \Z$.
Since the unique expansions which increase the number of $b$-edges involve the relation $cb^{-t}c\inv b^{t} = cb^{-t}b^{t}c\inv$, it follows easily by induction on the usual expansion/folding scheme that every edge in $X$ is $\sigma$-equivalent to some original edge. Thus $\sigma$ has only finitely many equivalence classes and we only need to show that the size of each equivalence class can be bounded. Let $k$ be the number of vertices in $\mt(u)$. Suppose that some equivalence class of $\sigma$ possesses $k+1$ different edges $p_i \xr{b} q_i$ for $i = 0,\ldots,k$. By (\ref{fpex4}), for each $i = 0,\ldots,k$ there
 is a path $q_i \xr{b^{m_i}} r_i$ in $S\Gamma(u)$ with $m_i \geq 0$ and $r_i$ original. Hence there exist $0 \leq i < j \leq k$ such that $r_i = r_j$. Since $q_i \neq q_j$, we get $m_i \neq m_j$. On the other hand, it is easy to see that, since $p_i \xr{b} q_i$ and $p_j \xr{b} q_j$ are $\sigma$-equivalent, there exists some path $q_i \xr{w} q_j$ in $S\Gamma(u)$ with $w\pi_b = 0$. But then
$$\xymatrix{
q_i \ar[rr]^{w} && q_j \ar[dl]^{b^{m_j}}\\
&r_i \ar[ul]^{b^{-m_i}}&
}$$
is a loop in $S\Gamma(u)$ with $(wb^{m_j}b^{-m_i})\pi_b \neq 0$, contradicting (\ref{fpex3}). Therefore each equivalence class of $\sigma$ has at most $k$ elements and so $S\Gamma(u)$ has finitely many $b$-edges.

It remains to bound the number of $c$-edges. Suppose that we remove all the $a$-edges and $b$-edges from $S\Gamma(u)$ to get the automaton $\A$. Each time we remove an edge, the number of connected components increases at most by one. Hence the number of connected components of $\A$ is bounded
and it suffices to show that all of them are finite.

Now  for every edge $p_i \xr{c} q_i \neq p_i$ appearing for the first time in $\A_i(u)$ in the Stephen's sequence, there exist necessarily some edge $p_{i-1} \xr{c} q_{i-1} \neq p_{i-1}$ in $\A_{i-1}(u)$ and paths
$$\xymatrix{
q_{i-1} && q_i \ar[ll]_{b^{t}}\\
p_{i-1} \ar[u]^c && p_i \ar[ll]^{b^{t}} \ar[u]_c
}$$
in $\A_i(u)$ because the expansion $ca = a$ could never produce the straight edge $p_i \xr{c} q_i$. It follows that, for every edge $p \xr{c} q \neq p$ in $S\Gamma(u)$, there exist some edge $p_1 \xr{c} q_{1}$ in $\mt(u)$ and paths
$$\xymatrix{
q_{1} && q \ar[ll]_{b^{m}}\\
p_{1} \ar[u]^c && p \ar[ll]^{b^{m}} \ar[u]_c
}$$
in $\A(u)$ for some $m \geq 0$. Since $\mt(u)$ is finite and $S\Gamma(u)$ has finitely many $b$-edges, then we can bound the size of each connected component of $\A$.
Thus $S\Gamma(u)$ possesses only finitely many $c$-edges and we are done.

(ii) For every $n \geq 1$, let $e_n = ((ab^{tn}a)(ab^{tn}a)^{-1})\p$.
We claim that
\beq
\label{fpex2}
\xymatrix{
p_0 \ar@(u,ur)^c \ar[r]_a & p_1 \ar@(u,ur)^c \ar[r]_{b^t} & p_2 \ar@(u,ur)^c \ar[r]_{b^t} &\ldots \ar[r]_{b^t}& p_{n} \ar@(u,ur)^c \ar[r]_{b^t} & p_{n+1} \ar@(u,ur)^c \ar[r]_a& p_{n+2}
}\eeq
is $S\Gamma(e_n)$. Here an edge $p\xrightarrow{b^t}q$  is shorthand for a sequence of $t$ edges in a straightline from $p$ to $q$, labeled by $b$, (together with their inverse edges) and  with no other edges incident on any vertex other than $p,q$.

It is clear that $\mt(e_n)$ is
$$\xymatrix{
\ar@{<->}[r] & p_0 \ar[r]_a & p_1 \ar[r]_{b^t} & p_2 \ar[r]_{b^t} & \ldots \ar[r]_{b^t}& p_{n+1} \ar[r]_a& p_{n+2}
}$$
Expanding through the relation $ca = a$ we can produce the first and the last $c$-loops. Then we expand through
the relation $cb^{-t}c\inv b^{t} = cb^{-t}b^{t}c\inv$ to successively produce the $c$-loops at $p_{n}, p_{n-1}, \ldots, p_1$. Since
(\ref{fpex2}) is an inverse automaton and admits no expansions, it must be $S\Gamma(e_n)$, which has therefore $tn+3$ vertices.

Now $\A(e_n)$ is obtained from $S\Gamma(e_n)$ by declaring $p_0$ the basepoint. Since $aca^{-1} \in L(\A(e_n))$, then $(aca^{-1})\p \geq e_n$ and so $e_n \in \lambda((aca\inv)\p) \cap E(S)$. We show that $e_n \in \mu_S((aca\inv)\p)$.

Suppose now that $f \in \lambda((aca\inv)\p) \cap E(S)$ satisfies $e_n \leq f$. Let
$v \in \wt{A}^*$ be a Dyck word representing $f$.

First we note that $f \in \lambda((aca\inv)\p)$ implies $f \leq (aca\inv)\p$. Hence $aca^{-1} \in L(\A(v))$.

Write $B = \{ b,c\}$. We show that

\beq
\label{fpex5}
\mbox{If $w \in \wt{B}^*$, then $S\Gamma(w)$ contains no $c$-loops.}
\eeq

First note that if $u\in \wt{B}^*$, then $u$ does not represent the same element of $S_t$ as any word $z$ containing $a$ or $a\inv$ since the relation $ca=a$ cannot be applied to any word in $\wt{B}^*$.  It then follows that the inverse subsemigroup $T_t$ of $S_t$ generated by $B$ is defined by the relation $cb^{-t}c\inv b^{t}=cb^{-t}b^{t}c\inv$.  It also follows that $a$ cannot label an edge in $S\Gamma(w)$ for any $w\in \wt{B}^*$  since if $a$ labels an edge in $S\Gamma(w)$, then $w\p = z\p$ for some word $z\in \wt{A}^+$ containing an $a$ or $a\inv$.  Thus the Sch\"utzenberger graphs of $w$ in $S_t$ and $T_t$ coincide.  From the presentation for $T_t$, $\pi_c|_{\wt{B}^*}$ factors through a homomorphism $T_t\to \mathbb Z$.  Since $\mathbb Z$ is a group, it follows that if $z\in \wt{B}^*$ labels a loop in $S\Gamma(w)$, then $z\pi_c=0$, and so (\ref{fpex5}) holds.

Since $e_n \leq f = v\p$, then $v \in L(\A(e_n))$. Let $\A$ be obtained by removing the vertex $p_{2n+2}$ from $\A(e_n)$ (in the notation given in (\ref{fpex2}) for $S\Gamma(e_n)$). Suppose that $v \in L(\A)$. We can factor
$$v = c^{i_0}(aw_1a\inv c^{i_1})\ldots (aw_ma\inv c^{i_m})$$
with $m\geq 1$, $i_j \in \Z$ and $w_j \in B^*$.  Expanding $\mathrm{Lin}(v)$ through the relation $ca = a$ and folding to get $c$-loops, we get the automaton
$$
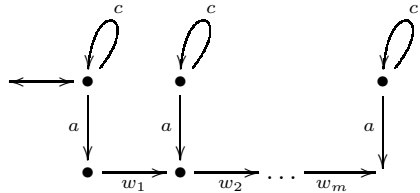
\xymatrix{
\ar@{<->}[r] & \bullet \ar@(ur,u)_c \ar[d]_a & \bullet \ar@(ur,u)_c \ar[d]_a&& \bullet \ar@(ur,u)_c \ar[d]_a \\
& \bullet \ar[r]_{w_1} & \bullet \ar[r]_{w_2} &\ldots \ar[r]_{w_m} &
}$$
Applying expansions and folds to the part corresponding to $\mathrm{Lin}(w_1\ldots w_m)$, the graph we get is $S\Gamma(w_1\ldots w_m)$ (which contains no $c$-loops by (\ref{fpex5})) with finitely many subgraphs
$\xymatrix{
p \ar@(u,ur)^c \ar[r]_a & q
}$
adjoined. No further expansion applies to this graph, hence what we get is really $\A(v)$. But we have remarked before that $aca^{-1} \in L(\A(v))$, contradicting the non existence of $c$-loops in $S\Gamma(w_1\ldots w_m)$. Thus
$v \in L(\A(e_n)) \setminus L(\A')$.

It follows that we can factor $v = v'v''$ so that $v'$ labels a path of the form
$$p_0\xr{v_0} p_0 \xr{a}  p_1 \xr{v_1} p_1 \xr{v_2} p_2 \xr{v_3} p_2 \xr{v_4} p_3 \xr{v_5} \ldots \xr{v_{2n}} p_{n+1} \xr{v_{2n+1}} p_{n+1}\xr{a} p_{n+2}$$
in $\A(e_n)$ for some factorization $v' = v_0av_1v_2\ldots v_{2n+1}a$ such that $v_{2j} \,\rho\, b^t$ for $j = 1,\ldots,n$. And we may assume that
\beq
\label{enh}
\mbox{the displayed edge $p_0 \xr{a}  p_1$ features the last occurrence of $p_0$ in the path labeled by $v'$.}
\eeq

For $j = 0,\ldots,n$, let $\A_j$ denote the inverse automaton depicted by
$$
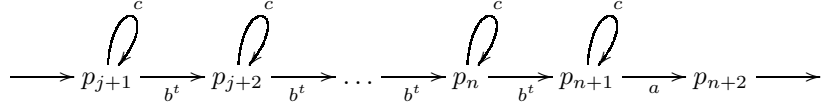
\xymatrix{
\ar[r] & p_{j+1} \ar@(u,ur)^c \ar[r]_{b^t} & p_{j+2} \ar@(u,ur)^c \ar[r]_{b^t} &\ldots \ar[r]_{b^t}& p_{n} \ar@(u,ur)^c \ar[r]_{b^t} & p_{n+1} \ar@(u,ur)^c \ar[r]_a& p_{n+2} \ar[r] &
}$$
We show that
\beq
\label{fpex8}
\mbox{$\A(w) = \A_{n-j}$ for all $j = 0,\ldots,n$ and $w \in L(\A_{n-j})$.}
\eeq

We use induction on $j$. The case $j = 0$ is immediate in view of the relation $ca = a$. Hence we assume that $j > 0$ and (\ref{fpex8}) holds for $j-1$.

Let $w \in L(\A_{n-j})$. Then we may write $w = w'xw''$ with $x\,\rho\,b^t$ and $w'' \in L(\A_{n-(j-1)})$. By the induction hypothesis, we get
$\A(w'') = \A_{n-(j-1)}$. Folding the $b$-edges and expanding through the relation $cb^{-t}c\inv b^{t} = cb^{-t}b^{t}c\inv$, we obtain $\A(b^tw'') = \A_{n-j}$.  But $w \in L(\A_{n-j})$ and $xw''\rho= b^tw''\rho$, thus $w\p\geq xw''\p$.  Therefore, $$w\p = (w'xw'')\p =(w'xw'')\p((xw'')\p)\inv (xw'')\p = w\p ((xw'')\inv xw'')\p =(xw'')\p$$
and so
$\A(w) = \A(b^tw'') = \A_{n-j}$, whence (\ref{fpex8}) holds.

In particular, $\A(z) = \A_0$ for $z = v_1b^tv_3b^t\ldots b^t v_{2n+1}a$ in view of (\ref{enh}). Expanding through the relation $ca = a$ and folding, it is easy to see that $\A(az)$ is just $\A(e_n)$ with $p_{2n+2}$ as terminal vertex. Since $v_0$ labels a loop at $p_0$, we have $v_0\p\geq (az)(az)\inv\p$ and so $v_0az\p = az\p$.  Thus $\A(v') = \A(v_0az) = \A(az)$. Hence $\A(v'(v')\inv) = \A(e_n)$ and so $(v'(v')\inv)\p = e_n$. But then
$$e_n = (v'(v')\inv)\p \geq v\p \geq e_n$$
yields $v\p = e_n$ and so $e_n \in \mu_S((aca\inv)\p)$ for every $n \geq 1$. Since we have established before that $S\Gamma(e_n)$ has $tn+3$ vertices, it follows that $\mu_S((aca\inv)\p)$ is infinite.
\medskip

\medskip

We can now use Example~\ref{fpex} to prove the following:

\bp
\label{imfp}
There exist infinitely many nonisomorphic finitely presented 3-generated inverse monoids satisfying the conditions of Problem~\ref{question}.
\ep

\proof
We have shown in Example~\ref{fpex} that $S_t$ satisfies the conditions of Problem~\ref{question} for every $t \geq 2$. Thus it suffices to show that $S_t \not\cong S_n$ for distinct $t,n \geq 2$.

Suppose that $\theta\colon S_t \to S_n$ is an isomorphism. It is easy to see that each generating set of $S_n$ must contain
\bi
\item
$a^{\alpha}$ with $\alpha = \pm 1$;
\item
$b^{\varepsilon}$ with $\varepsilon = \pm 1$;
\item
$c^{\delta}$ with $\delta = \pm 1$.
\ei
Hence any minimal generating set of $S_n$ consists necessarily of three elements of this form.Thus we may assume that $\theta$ is induced by some bijection $\{ a,b,c\} \to \{ a^{\alpha},b^{\varepsilon},c^{\delta} \}$.

Since $ca = a$ is a relation of the presentation of $S_t$, then
$(c\theta)(a\theta) = a\theta$ holds in $S_n$, yielding successively $c\theta = c^{\delta}$ (since $a$ and $b$ never label loops), $a\theta = a^{\alpha}$ and $b\theta = b^{\varepsilon}$.

But $cb^{-t}c\inv b^{t} = cb^{-t}b^{t}c$ is also a relation of the presentation of $S_t$, hence
$c^{\delta}b^{-\varepsilon t}c^{-\delta}b^{\varepsilon t} = c^{\delta}b^{-\varepsilon t}b^{\varepsilon t}c^{-\delta}$ holds in $S_n$. Since $\mt(c^{\delta}b^{-\varepsilon t}c^{-\delta}b^{\varepsilon t})$ admits no expansion for the presentation of $S_n$,
we have reached a contradiction. Therefore $S_t \not\cong S_n$.
\qed

\section*{Acknowledgments}

The first author acknowledges support from the Center of Mathematics of the University of Porto, which is financed by national funds through the Funda\c c\~ao para a Ci\^encia e a Tecnologia, I.P., under the project with reference UIDB/00144/2020. The second author was supported by a PSC CUNY grant and a Simons Foundation Collaboration Grant, award number 849561.

\end{document}